\documentstyle[12pt]{article}
\renewcommand{\theequation}{{\thesection}.{\arabic{equation}}}
\def\r#1{\mbox{(}\ref{#1}\mbox{)}}
\def\ot{\otimes}
\def\kEi{k(E_{\la_i})}
\def\Ei{E_{\la_i}}

\def\kMl{k(\Ml)}

\def\RVM{{\cal R}{(V,\, {M_{\mu}^q)}}}
\def\RVVM{{\cal R}{(V^{\ot 2},\, {M_{\mu}^q)}}}
\def\RVMh{{\cal R}{(V,\, {M_{\mu}^{q\h})}}}
\def\kMm{k({\overline M}_{\mu})}
\def\RVm{{\cal R}(V, {\overline M}_{\mu})}
\def\Rn{{\cal R}_{\nu}}
\def\Ru{{\cal R}_{\nu}}
\def\Rmi{{\cal R}_{\mu_i}}

\def\la{\mu}
\def\Ml{{M_{\la}}}

\def\uq{U_q({\bf g})}
\def\uqs{U_q(sl(n))}

\def\Sym{\rm Sym}

\def\h{{\hbar}}

\def\VV{V^{\ot 2}}

\def\h{{\hbar}}

\def\bea{\begin{eqnarray}}
\def\eea{\end{eqnarray}}
\def\beq{\begin{equation}}          
\def\eeq{\end{equation}}

\def\rank{{\rm rank\, }}

\def\Im{{\rm Im\, }}
\def\Id{{\rm Id }}

\def\Id{{\rm id\, }}
\def\Vect{{\rm Vect\,}}

\def\Tr{{\rm Tr}}
\def\det{{\rm det}}

\def\dim{{\rm dim}\,}
\def\id{{\rm id}\,}
\def\span{{\rm span}\,}

\def\al{\alpha}
\renewcommand{\theequation}{{\thesection}.{\arabic{equation}}}
\newtheorem{proposition}{Proposition}

\newtheorem{theorem}{Theorem}

\newtheorem{remark}{Remark}
\def\R{{\cal R}}
   
 \def\g{{\bf g}}
\def\pole{k}
\def\L{{\cal L}(P)}
\def\Lq{{\cal L}_q(R)}
\def\Lh{U(gl(n)_{\hbar})}
\def\Lqh{{\cal L}_{q,\hbar}(R)}
\def\Lhh{{\cal L}_{\hbar}}
\def\Tq{{\cal T}_q(R)}

\title{Quantum line bundles\\
via Cayley-Hamilton identity}
\author{D.~Gurevich,\\
{\footnotesize\it ISTV, Universit\'e de Valenciennes,
59304 Valenciennes, France}\\
P.~Saponov\\
{\footnotesize\it Theory Division, Institute for High Energy
Physics, 142284 Protvino, Russia}}
\begin{document}
\maketitle
\begin{abstract}
As was shown in \cite{GPS} the matrix $L=\vert\vert l_i^j\vert\vert$
whose entries $l_i^j$ are generators of the so-called reflection
equation algebra is subject to some polynomial identity looking like
the  Cayley-Hamilton identity for a numerical matrix. Here a similar
statement is presented for a matrix whose entries are generators of a
filtered algebra being a "non-commutative analogue" of the reflection
equation algebra. In an appropriate limit we get a similar statement
for  the matrix formed by the generators of the algebra $U(gl(n))$.
This property is used to introduce  the notion of line bundles over
quantum orbits in the spirit of the Serre-Swan approach. The quantum
orbits in question are presented explicitly as some quotients of one
of the mentioned above algebras both in the quasiclassical case (i.e.
that related to the quantum group $U_q(sl(n))$) and a
non-quasiclassical one (i.e. that arising from a Hecke symmetry with
non-standard Poincar\'e series of the corresponding symmetric and
skewsymmetric algebras).
\end{abstract}



\section{Introduction}

 Let $G=GL(n),\,\,{\g}=Lie\,(G)$ over the field $\pole={\bf R}$ or
 $\pole ={\bf C}$\footnote{The choice of the field is similar to that
 in the classical case. If $k={\bf R}$ the entries of the quantum
 R-matrix $R$ are assumed to be real.}. As usual, let us identify
 ${\g}$ and ${\g}^*$ and consider a matrix $A\in {\g}^*$ with pairwise
 distinct eigenvalues $\la_1,\,\la_2,\,...\la_n$. Denote
 $\Ml,\,\,\la=(\la_1,\,\la_2,\,...\la_n)$ a $G$-orbit of the matrix
 $A$ w.r.t.  coadjoint action of $\g$. Then the orbit $\Ml$ being an
 affine algebraic variety is defined by the following system of
 equations
 \beq
 {\rm \Tr}\, A^k=c_k=\sum \la_i^k,\,\,\, k=1, 2, 3,...,n.
 \label{odin}
 \eeq
  Let $\kMl$ be its coordinate ring.

 A typical example of a line bundle over $\Ml$ is an eigenspace
 corresponding to an eigenvalue $\la_i$, i.e. the space of the vectors
 $v\in V^*$ such that
 \beq
 v\,A= \la_i\, v,\,\, v\in V^*\;,
 \label{dva}
 \eeq
 $V^*$ being the right $\g^*$-module.

 Then this line bundle itself is an algebraic variety: it is defined
 in the space $\g^*\times V^*$ by system \r{odin}---\r{dva}. This
 variety (i.e. the total space of the line bundle in question) will be
 denoted $\Ei$. The coordinate ring of this variety $\kEi$ has the
 structure of a $\kMl$-module.

 This example is a particular case of a one-to-one correspondence
 between  algebraic vector bundle over an affine algebraic variety and
 finitely generated projective modules over its coordinate ring
 realized in \cite{Se} (a similar correspondence on compact smooth
 varieties was established in \cite{Sw}).

 Let us remark that  the projectivity of the $\kMl$-module $\kEi$ can
 be shown by means of the projector \beq P_i = \prod_{j\not=
 i}^n{{A-\la_j\,\id}\over{\la_i-\la_j}} \label{proj} \eeq since this
 $\kMl$-module can be identified with $\Im\,P_i$.

 Our main purpose is to generalize the construction of considered (and
 some "derived") line bundles to the quantum case.

 Up to our knowledge the first attempt to realize a line bundle over a
 quantum sphere  in terms of projective modules featured in \cite{HM}.
 Constructed there was a quantum  analogue of projector \r{proj} over
 a quantum sphere (or what is the same,  a quantum hyperboloid if we
 ignore the involution operation).

 In this paper we suggest  a regular way of constructing projective
 modules for a  generic  $\la$, which are the quantum deformations of
 $\kMl$-modules $\kEi$. More precisely, we construct quantum two
 parameter deformation of the algebra $\kMl$ and  line bundles $\kEi$
 (see below).

 The basic question  arising from the very beginning is what are
 quantum analogues of the orbits in $\g^*$. An habitual way to
 introduce such quantum objects makes use of so-called Hopf-Galois
 extension (cf. \cite{Sh}, \cite{HM} ). This approach allows one to
 generalize the notion of an ordinary orbit in terms of a couple of
 Hopf algebras. The famous "RTT" algebra  and some its Hopf subalgebra
 are usually employed as such couples. However, such an approach does
 not enables one either to control the flatness of
 deformation\footnote{We refer the reader to \cite{DGK} for the
 rigorous definition of this notion. Roughly speaking, this means that
 the supply of elements does not change under deformation of the
 initial object. Let us remark that under a flat deformation $A_h$ of
 a
 commutative associative algebra $A=A_0$ the linear in $h$ 
skewsymmetrized term of
 the deformed product is a Poisson bracket.} in the quasiclassical
 case or
 to generalize construction of quantum orbits to a non-quasiclassical
 case.

 Let us note that we use the term "quasiclassical" for objects arising
 from deformations of classical ones. In this paper all quasiclassical
 objects in question are endowed with a structure of $\uqs$-module
 (and as usual the products in all quantum algebras involved are
 supposed to be $\uqs$-covariant). By "non-quasiclassical" objects we
 mean those arising from the solutions $R$ of the quantum Yang-Baxter
 (YB) equation \r{YB} whose symmetric and skewsymmetric algebras have
 non-quasiclassical Poincar\'e series. These algebras are well defined
 if we assume $R$ to be a Hecke symmetry (i.e. a solution of the YB
 equation subject to the Hecke condition \r{Hecke}). For details the
 reader is referred to \cite{G} where a large family of
 non-quasiclassical Hecke symmetries was constructed.

We suggest another way of introducing quantum orbits. The central
role in our approach is played by so-called reflection equation (RE)
algebra $\Lq$ (see Section 2 for definition) which can be associated
with any (quasiclassical or not) Hecke symmetry $R$. If $R$ is a
quasiclassical Hecke symmetry then the algebra $\Lq$ (similarly to the
RTT algebra) is a flat deformation of the coordinate ring
$k(Mat(n,\,k))$. However, its properties differ drastically from those
of the RTT algebra.

 The main difference is that the RE algebra possesses a big center
 $Z$. In particular, the quantum trace $\Tr_q$ well defined in this
 algebra belongs to $Z$. On  quotienting the RE algebra over the ideal
 generated by $\Tr_q$ we get an algebra with $n^2-1$ generators which
 we consider as  a q-analogue of the algebra $k(sl(n)^*) = {\Sym}
 (sl(n))$. If instead of the mentioned ideal we take that generated by
 the elements
 \beq
 z-\chi(z),\,\, z\in Z \label{center}
 \eeq
 where $\chi$ is a generic character of $Z$, we get the quotient
 algebra (denoted $k(M_{\mu}^q)$) which can be considered as a quantum
 analogue of semisimple orbits above (we call semisimple orbits those
 of semisimple elements). Thus, both in the quasiclassical and
 non-quasiclassical cases such type quantum orbits are defined by
 means of some "quantum (or braided) algebraic equations"  in the
 spirit of affine algebraic geometry. By abusing the language we call
 quantum orbits the corresponding "coordinate rings" in both cases.

 Another important difference between the RTT and RE algebras consists
 in the fact that the latter one (being a quadratic algebra) admits a
 further flat deformation giving rise to a quadratic-linear algebra
 which looks like the enveloping algebra $U(gl(n))$. The final object
 of such a deformation is an algebra $\Lqh$ (see Section \ref{deform})
 depending on two parameters which tends to the RE algebra as $\h\to
 0$ and to $U(gl(n)_{\h})$ as $q\to 1$, where the defining relations
 of $U(gl(n)_{\h})$ are as follows:
 \beq
 a^{i_1}_{j_1}a^{i_2}_{j_2} - a^{i_2}_{j_2}a^{i_1}_{j_1} =
 \hbar(\delta^{i_2}_{j_1}a^{i_1}_{j_2}-\delta^{i_1}_{j_2}
 a^{i_2}_{j_1})\;.\label{commut}
 \eeq
 Hereafter we use the notation $U(\g_{\h})$ for the enveloping algebra
 of a Lie algebra $\g_{\h}$ which differs from $\g$ by the factor $\h$
 introduced into the Lie bracket. We prefer to use the Lie algebra
 $\g_{\h}$ instead of $\g$ in order to represent its enveloping
 algebra $U(\g_{\h})$ as a deformation of the algebra
 $k(\g^*)=\Sym(\g)$.  In a similar way we treat quotients of the
 algebra $U(\g_{\h})$ as deformations of the corresponding
 orbits\footnote{As for the RTT algebra, it does not have any
 non-trivial quadratic-linear deformation which could be considered 
 as a q-analogue of $U(gl(n))$ (cf. \cite{GR}).}.

 Similarly to $\Lq$ the algebra $\Lqh$ has a big center. On
 quotienting the algebra $\Lqh$ over an ideal looking like \r{center}
 we get "quantum non-commutative" analogue of the orbits above. We
 treat the specialization of this quotient at the point $q=1$ as a
 "classical non-commutative" orbit. Thus, this specialization is just
 an appropriate quotient of the algebra $U(gl(n)_{\h})$ (or
 $U(sl(n)_{\h})$).

 Note that we consider the RE algebra $\Lq$ and its quotients as "
 quantum commutative" algebras. Their non-commutative counterparts are
 the algebra $\Lqh$ and its quotients (all these algebras are well
 defined in the non-quasiclassical case as well). The term "classical"
 means that the product in the algebra in question is $G$-covariant
 where $G$ is a usual group. By contrast, "quantum" means that the
 product in the algebra in question is covariant w.r.t. a Hopf
 algebra. In the quasiclassical case this Hopf algebra is just $\uqs$.
 In a non-quasiclassical case an explicit description of a similar
 Hopf algebra is more complicated (cf. \cite{AG} where an attempt to
 describe such an algebra featured). This is the reason why in a
 generic case it is more convenient to use the RTT algebra in order to
 define "symmetries" of the objects in terms of its coaction (see
 Section  \ref{deform}).

Now  let us explain what we understand by quantum analogues of the
line bundles above. Replace the matrix $A$ in \r{dva}  by matrix
$L=\vert\vert l_i^j\vert\vert$ subject to \r{REA} (this means that the
matrix $L$ is formed by the elements $l_i^j$ satisfying the quadratic
relations \r{REA}). Thus, we have the following system
 \beq
 v_i\,l^i_j-\nu\, v_j=0,\qquad \nu\in k, \label{qb}
 \eeq
where the summation over repeated indices is assumed. Otherwise
stated, we consider the free right $\Lq$-module
$$
V\ot_k\Lq\,\,\,\,{\rm where}\,\,\,\, V=\span (v_i)
$$
and its submodule ${\cal R}_{\nu}$ generated by the  l.h.s. of \r{qb}.
Let us restrict ourselves to a "quantum orbit" $k(M_{\mu}^q)$. This
means that instead of the free $\Lq$-module $V\ot_k\Lq$ and its
submodule $\R_{\nu}$ we consider the free right $k(M^q_{\mu})$-module
$$\RVM=V\ot k(M^q_{\mu})$$ and its submodule generated by the l.h.s.
of \r{qb} (we keep the notation $\R_{\nu}$ for it).

We call {\em a quantum line bundle} over the given quantum orbit
$k(M_{\mu}^q)$ the quotient
$$
\RVM/{\cal R}_{\nu}
$$
if it is non-trivial. This definition can be extended to
quantum line bundle over "non-commutative quantum orbits". For
that its suffices to replace "commutative quantum orbit" in this
definition by its "non-commutative" counterpart (denoted
$k(M_{\mu}^{q\h})$).  The problem is when the quotient
$\RVM/{\cal R}_{\nu}$ (or its non-commutative analogue
$\RVMh/{\cal R}_{\nu}$) is non-trivial.  In the classical
commutative case $(q=1,\,\h=0)$ it is so iff $\nu$ is a root of
the characteristic polynomial of $A$ (i.e., $\nu=\mu_i$ for some
$i$).

In this paper we give a criterion on $\nu$ which yields
non-triviality of these quotients.  This criterion is based on
the quantum version of the Cayley-Hamilton (CH) identity for the
matrix $L$ found in \cite{GPS}. This version of the CH theorem
states that there exists a polynomial $P$ whose coefficients
belongs to $Z$ and such that $P(L)=0$. When we restrict
ourselves to a quantum orbit the  coefficients of $P$ become
numerical. So, we get a polynomial $\overline P$ with numerical
coefficients such that ${\overline P}(L)=0$ . Our main statement
says that the quotient module $\RVM/{\cal R}_{\nu}$ is
non-trivial iff $\nu$ is a root of $\overline P$ (we assume that
the roots $\mu_i$ of $\overline P$ are pairwise distinct).
Moreover, this quotient is projective and  in the quasiclassical
case it is a flat deformation of its classical counterpart.

Besides, we present here a version of the CH identity valid for
two parameter algebra $\Lqh$ and by passing to the limit $q\to
1$ we get such an identity for the algebra $U(gl(n)_{\hbar})$.
This allows us to  get a similar description for
"non-commutative orbits" both in the classical and quantum
cases.  Let us remark that a version of the CH identity for the
algebra $U(gl(n))$ is already known from the late sixties due to
works \cite{BiLo}. But in the cited works the identity was
established for any finite representation of $U(gl(n))$ and
looks like
$$
\prod_i (A-\mu_i) = 0\;,
$$
where $\mu_i$ are {\em integer} numbers, depending on
 a given representation. In a sense, our result is
more general, since the CH identity is realized with
coefficients being elements of  the RE algebra itself without
using any representation. In particular, this allows us to
consider the orbits of general form, where the coefficients
$\mu_i$ are not obligatory integer numbers.
Also a non-commutative version of the CH identity is presented
in \cite{GKLLRT}.  However it is rather useless for our aims
since the coefficients of the CH  polynomial are scalar
matrices.

In the classical case besides the above line bundles related to
the fundamental vector $sl(n)$-module  $V$ (called in the sequel
{\em basic}) there exist other line bundles which can be
obtained via the tensor products of the basic ones.  Moreover,
the family of all line bundles over a regular algebraic (or a
smooth) variety forms a ring w.r.t.  the tensor product.  Then a
natural question arises: what is a regular way to construct
quantum line bundles over the "quantum orbits" which would be
different from the basic ones  (we will refer to them as {\em
derived} line bundles).  If we want to realize the tensor
product of two or more basic line bundles in terms of projective
modules we should construct the corresponding projector.
In fact the problem of constructing such a projector reduces to
the problem of finding the CH identity for the matrix $L$
extended to the tensor product of two (or more) copies of the
space $V$. In the classical case this CH identity can be easily
found. 

However, it is not so  in the quantum case. It is not
even clear what is a reasonable way to extend the action of the
matrix $L$ to a tensor power of the space $V$.  Remark that any
"reasonable" at first glance way leads to an extension of the
matrix $L$ for which we are not able to find any polynomial
identity which would be a deformation of the classical one (see
Section 5).  Nevertheless, there exists a "canonical" way to
extend the action of the matrix $L$ to the {\it symmetric\/}
part of $V^{\ot l}$.  Hopefully,  for such an extension of the
matrix $L$  the CH identity can be found and it is a flat
deformation of its classical counterpart. At least, it is so in
a particular case when $l=2$ and $\rank (R)=2$ (this means that
the Poincar\'e series of "skewsymmetric algebra" of the space
$V$ is of the form $P_-(t)=1+n\,t+ t^2$).

In subsequent  publications we will apply our approach to a
quantum version of K-theory which on one hand would enable us to
control the flatness of deformation in the quasiclassical case
and on the other hand would be valid in non-quasiclassical case.

The paper is organized as follows. In Section 2 we give a
description of the RE algebra in  comparison with RTT algebra
and introduce "quantum orbits" as some quotients of the former
algebra.  "Non-commutative" counterparts of these orbits are
introduced as well.  In Section 3 we present quasiclassical
counterparts of these quantum orbits assuming them to be
deformations of generic semisimple ordinary orbits. Section 4 is
devoted to construction of "basic line bundles" over quantum
orbits in terms of projective modules. And in the last Section
we discuss a way to define some derived line bundles related to
the symmetric product of the basic modules and calculate the
corresponding CH identity in the simplest case mentioned above.
In Appendix we present calculations of coefficients of the CH
identity for "non-commutative" cases.

{\bf Acknowledgment} The authors are supported by the grant
PICS-608/RFBR 98-01-22033. One of the authors (P.S.) is recognized to
Laboratoire de Math\'ematiques de l'Universit\'e de Valenciennes
for hospitality, during his stay at this laboratory the paper
was started.

\section{Reflection equation algebra and quantum orbits}
\label{deform}
\setcounter{equation}{0}

Consider a matrix solution $R^{i_1i_2}_{j_1j_2} \in {\rm
Mat}(n^2,\pole)$ of the Yang-Baxter equation\footnote{The standard
matrix conventions of \cite{FRT} are used throughout the paper.}
\beq
R_{12}R_{23}R_{12} = R_{23}R_{12}R_{23}\;,
\label{YB}
\eeq
satisfying the additional Hecke condition
\beq
R^2 = {\rm id } + \lambda R\qquad {\rm where} \quad \lambda =
q - q^{-1}\;,
\label{Hecke}
\eeq
the value of nonzero number $q\in \pole$ being generic:
$q^r\not=1$ for any integer $r$. Such solutions will be refered
to as Hecke symmetries, and following to \cite{G} we will also
suppose, that the Hecke symmetry $R$ is an even symmetry of
finite rank $p\le n$. This means that
$$
P_-^{(p+1)}\equiv 0\quad {\rm and}\quad {\rm dim}\,P_-^{(p)}=1\;,
$$
where   $P^{(l)}_{-}$ stands for the projector of $V^{\ot l}$
onto its subspace of totally skewsymmetric tensors.  It is
possible to show that such a Hecke symmetry is closed in the
sense of \cite{G}, i.e.,  the matrix $R^{t_1}$ is invertible.
More detailed treatment can be found in \cite{G,GPS}.

With any Hecke symmetry $R$ (quasiclassical or not) we can
associate two matrix algebras. One of them, denoted below as
$\Tq$ and called RTT algebra is generated by $n^2$ quantities
$t^i_j$ which can be considered as entries of some matrix $T
=\vert \vert t^i_j\vert\vert$ subject to the following quadratic
relations \cite{FRT}:
\beq
R_{12} T_1 T_2 = T_1 T_2 R_{12}\;.\label{RTT}
\eeq
It is well known that the algebra $\Tq$  possesses the
bialgebra structure w.r.t. to the comultiplication:
$$
\Delta(t^i_j) = t^i_p\otimes t^p_j\;.
$$
In the case of an even Hecke symmetry (quasiclassical or not)
and under the assumption that so-called quantum determinant is
central (cf. \cite{G}) we can extend this bialgebra structure to
the Hopf algebra by introducing the antipodal mapping
$$
{\cal S}:\ \Tq\to\Tq.
$$
Another of mentioned algebras is so-called RE algebra
$\Lq$\footnote{There are known different versions of the RE
algebra, cf. \cite{KulSkl,KulSas}.  We use that introduced in
\cite{M1} in terms of braided matrix algebra.} The corresponding
$n\times n$ matrix $L=||l^i_j||$ obeys the relation
\beq
R_{12}L_1R_{12}L_1 = L_1R_{12}L_1R_{12}\;.\label{REA}
\eeq

This algebra can be given a structure of the adjoint comodule
w.r.t. the coaction $\delta$ of the algebra $\Tq$:
\beq
\delta(l^{i}_j) = t^i_p \,{\cal S}(t^k_j)\otimes l^p_k.\label{comod}
\eeq
\begin{remark} As we have said in Introduction we prefer using
the RTT algebra as a substitute of the symmetry group since it
is well defined in the both quasiclassical and
non-quasiclassical cases. As for the dual object, its explicit
description in a non-quasiclassical case is not easy (cf.
\cite{AG}). However, in fact we can do without the RTT algebra
at all.  Let us also point out that the RE algebra has the
structure of a braided Hopf  algebra. This property was
discovered by S.Majid (cf. \cite{M2}). However, we do not use
this property either.
\end{remark}

An important feature of the both algebras mentioned above is the
existence of polynomial identities on the quantum matrices  $L$
and $T$ \cite{GPS,IOP}, which generalize the well-known CH
identity of the classical matrix analysis.

For the RE algebra $\Lq$ this identity looks as follows \cite{GPS}:
\beq
(-L)^p +\sum_{k=0}^{p-1} (-L)^k\sigma_{p-k}(L) \equiv 0\;.
\label{RE-CH}
\eeq
Note, that the upper limit of summation in the identity is defined by
the number $p={\rm rank}(R)$, not by $n$ --- the dimension of the
space $V$.  The coefficients $\sigma_k(L)$ in the above relation are
polynomial combinations of generators $l^{i}_j$:
\beq
\begin{array}{l}
\sigma_k(L) = \alpha_k\, {\rm Tr}_{(12\dots p)}
P_-^{(p)}(L_1R_{12}\dots R_{k-1,k})^k\\
\\
\alpha_k = q^{-k(p-k)} [C^k_p]_q
\end{array}
\label{L-sig}
\eeq
where $[C^k_p]_q$ are $q$-binomial coefficients. Let us note that in
the quasiclassical case the set $\{\sigma_k(L)\}$ generates the center
$Z$ of the  algebra  $\Lq$. We will conjecturally suppose the same
property to be true in the non-quasiclassical case as well.

The CH identity for the algebra $\Tq$ is completely different from
that above.  As was shown in \cite{IOP} the matrix $T$ of generators
of  $\Tq$ algebra satisfies the following identity:
\beq
(-T)^{\bar p} + \sum_{k=1}^{p-1}(-T)^{\bar k}\sigma_{p-k}(T)
+ \sigma_{p}(T){\cal D} \equiv 0\;,\label{RT-CH}
\eeq
where ${\cal D}$ is a numeric matrix and
\bea
&&T^{\bar k} = {\rm Tr}_{(2\dots k)}R_{12}R_{23}\dots
R_{k-1,k}T_1T_2\dots T_k\\
&&\sigma_k(T) = \alpha_k\, {\rm Tr}_{(12\dots
p)}P_-^{(p)}T_1T_2\dots
T_k\;.\label{T-sig}
\eea

In contrast to the  algebra  $\Lq$, the  quantities
(\ref{T-sig}) are not central, they only form a commutative
subalgebra of $\Tq$.  It is this property that prevents us from
defining a quantum orbit in $\Tq$ algebra as a quotient algebra
over an ideal generated by the elements $\sigma_k(T)-c_k$ since
due to the non-centrality of $\sigma_k(T)$ the corresponding
quotient would not be a flat deformation of its classical
counterpart.

Now we consider a special case of the quasiclassical Hecke
symmetry related to the QG $\uqs$ (see Introduction). In this
case the $R$-matrix is a deformation of the  usual permutation
$P:\ P_{12}(v_1\otimes v_2) = v_2\otimes v_1$.  This means that
$\lim_{q\rightarrow 1}R = P$.  Therefore, at the limit
$q\rightarrow 1$ the quadratic quantum algebra (\ref{REA}) turns
into  the commutative algebra $\L = \lim_{q\rightarrow 1}\Lq$:
$$
P_{12}L_1P_{12}L_1 - L_1P_{12}L_1P_{12} \equiv L_1L_2 - L_2L_1 =
0\;
$$
where as usual, $L_2=P\,L_1\, P$.
Let us pass from the quadratic algebra $\Lq$ to a {\it
quadratic-linear\/} algebra $\Lqh$ with two parameters  $q$ and
$\hbar$ which can be treated as deformation of $\Lq$. Note, that
this deformation is also well defined in a non-quasiclassical
case. The algebra $\Lqh$ will be introduced by the following
simple procedure.  On shifting the generators of $\Lq$  $l^{i}_j
= {\bar l}^{i}_j - h \delta^{i}_j$ we come to the equivalent
algebra:
$$
R_{12}{\bar L}_1R_{12}{\bar L}_1 - {\bar L}_1R_{12}{\bar
L}_1R_{12} = \lambda h (R_{12}{\bar L}_1 - {\bar L}_1R_{12})\;.
$$
Now redefining the combination $\lambda h$ as a new parameter
$\hbar$ and treating it as {\it independent \/} on $q$, we get
the two-parameter quadratic-linear algebra  $\Lqh$
\beq
R_{12}{\bar L}_1R_{12}{\bar L}_1 - {\bar L}_1R_{12}{\bar
L}_1R_{12} = \hbar (R_{12}{\bar L}_1 - {\bar L}_1R_{12})\;.
\label{Lqh-com}
\eeq
We keep the notation $\bar L$ for the quantum matrix formed by
the generators of algebra $\Lqh$  in order to distinguish it
from the matrix $L$.
In the quasiclassical case the algebra $\Lqh$ can be considered
as a two parameter deformation of the commutative algebra $\L=
k(gl(n)^*) = \Sym(gl(n))$ since as is evident from
(\ref{Lqh-com})
\beq
\L = \lim_{q\rightarrow 1 \atop \hbar\rightarrow 0} \Lqh \qquad
\Lq = \lim_{\hbar\rightarrow 0 \atop q={\rm const}} \Lqh
\;.\label{limits}
\eeq

It is important, that all these deformations are flat. Otherwise
stated, the Poincar\'e series of $\L$, $\Lq$ and that of the
graded algebra associated to $\Lqh$ are equal to each other.
This statement is also valid in the non-quasiclassical case if
limits \r{limits} exist. Besides, the algebra $\Lqh$ admits a
nontrivial classical limit --- the noncommutative algebra
\beq
\Lhh=\Lh=\lim_{q\rightarrow 1 \atop \h={\rm const}} \Lqh.
\label{h-lim}
\eeq
Indeed, as $q\to 1$ the quadratic-linear relations (\ref{Lqh-com})
turns  into the following ones
\beq
A_1A_2 - A_2A_1 = \hbar(A_1P_{12}-P_{12}A_1)\;,\label{U-com}
\eeq
which are just relations \r{commut}. Moreover, we have
$$
\L=\lim_{\hbar  \rightarrow 0} \Lhh\;.
$$

Accordingly to what we said in Introduction in the
quasiclassical case we treat the algebras $\Lq$,$\Lqh$, and
$\Lh$ respectively as "quantum commutative", "quantum
non-commutative" and "classical non-commutative" counterparts of
the commutative algebra $\L$.

The crucial point is that for the algebras $\Lqh$ and $\Lh$
there also exists some version of the CH identity.  For the
matrix $\bar L$ formed by the generators of the algebra $\Lqh$
the CH identity is similar to \r{RE-CH}:
\beq
(-{\bar L})^p + \sum_{k=0}^{p-1} (-{\bar L})^k
\sigma^{(\hbar)}_{p-k}({\bar L}) \equiv 0\label{Lqh-CH}
\eeq
(this relation is valid in a non-quasiclassical case as well).
By passing to the limit $q\to 1$ we get a version of the CH
identity for the matrix $A$ formed by the generators of $\Lh$:
\beq
(- A)^p + \sum_{k=0}^{p-1}(-A)^k\tau^{(\hbar)}_{p-k}( A) \equiv
0\label{Ugl-CH}
\eeq
(see  Introduction on other versions of the CH identity).
The coefficients $\sigma^{(\hbar)}_k$ and $\tau^{(\hbar)}_k$ are
central elements of the corresponding algebras. The explicit
form of coefficients $\sigma^{(\hbar)}_k$ and the proof of the
existence of their limits (denoted $\tau^{(\hbar)}_k$) as $q\to
1$ are presented in Appendix.

Let us pass now to quantum analogues of generic orbits in
question.  We introduce such a "quantum orbit" (both in the
quasiclassical and non-quasiclassical cases) as the  quotient of
the RE algebra over the ideal generated by the elements
\r{center}. Keeping in mind the fact that the center $Z$ is
generated by the elements $\sigma_k(L)$ we can define the character
$\chi$ by imposing $\chi(\sigma_k(L))=c_k$. In other words, we define
the quantum orbit in question by the system of polynomial equations:
\beq
\sigma_k(L) - c_k=0\qquad k=1\dots p\;.\label{sig-fix}
\eeq
Then the  polynomial $\overline P$ mentioned in Introduction becomes
\beq
{\overline P}=(-L)^p +\sum_{k=1}^p(-L)^k c_{p-k} \;. \label{poly}
\eeq

We consider the quotient of the RE algebra over the ideal
generated by the l.h.s. of \r{sig-fix} as a quantum  commutative
orbit.  Let us denote this quotient by $k(M_{\mu}^q)$,
$\mu=(\mu_1,...\mu_p)$, $\mu_i$ being roots of the polynomial
$\overline P$.  Let us remark that in the quasiclassical case
the quotient $k(M_{\mu}^q)$ is a flat deformation  of its
classical counterpart $k(M_{\mu})$ being the coordinate ring of
a semisimple orbit. It can be shown by the  methods of
\cite{D1}.  In a similar way we can introduce non-commutative
quantum orbit $k(M_{\mu}^{q\h})$ as the quotient of the algebra
$\Lqh$ over the ideal generated by $\sigma_k^{(\h)}-c_k$ (with a
similar meaning of $\mu$).  The corresponding  polynomial
\r{poly} can be obtained if we replace the coefficients
$\sigma^{(\hbar)}_{k}({\bar L})$ in \r{Lqh-CH} by $c_k$.  In the
quasiclassical case by passing to the limit $q\to 1$ we get a
classical non-commutative orbit (its "coordinate ring" will be
denoted $k(M_{\mu}^{\h})$) and the corresponding polynomial
\r{poly}.  In all cases we suppose the roots $\mu_1,...\mu_n$ of
corresponding polynomials to be  pairwise distinct.

\section{Quasiclassical case: related  Poisson structures}
\setcounter{equation}{0}
In this section  we will briefly describe quasiclassical
counterparts of the algebras $\Lq$ and $\Lqh$ and their
restrictions to the orbits in question assuming $R$ to be
quasiclassical Hecke symmetry (see Introduction).  Let
\beq
r=\sum X_{\al}\wedge X_{-\al}\in \wedge^2(\g)
\label{rmat}
\eeq
be the classical r-matrix related to a simple classical Lie
algebra $\g$. The quasiclassical counterpart of the RTT algebra
is well known.  It is so-called Sklyanin bracket defined as the
difference between the left-invariant and  right-invariant
brackets on the corresponding Lie group $G$  associated to
r-matrix \r{rmat}, these invariant brackets are not Poisson
separately. Sklyanin bracket can be reduced to any semisimple
orbit  in $\g^*$. It can be quantized in the sense of
deformation quantization (cf. \cite{DG}). The resulting algebra
is $\uq$-covariant.  This quantum algebra can be also treated in
terms of the Hopf-Galois extension mentioned in Introduction.

However, the reduced Sklyanin bracket well defined on any
semisimple orbit is not defined  on the whole of $\g^*$. Roughly
speaking, we say that the Poisson brackets defined on each
semisimple orbit separately cannot be "glued"  into a global
Poisson bracket.

By contrast, the quasiclassical counterpart of the RE algebra is
well defined on the whole $gl(n)^*$. Let us describe it.  Let us
put $\g=sl(n)$ and associate to r-matrix \r{rmat} a bi-vector
field arising from the representation
$$
ad^*:\g\to \Vect(\g^*).\
$$
By applying this bi-vector field to functions $f$ and $g$ we get
a bracket $\{f,\, g\}_r$ which is not Poisson. Nevertheless, by
adding  some invariant summand we can convert it into a  Poisson
bracket. This summand can be constructed as follows.  It is well
known that in the decomposition of $\g^{\ot 2}$ into a direct
sum of irreducible $\g$-modules the component isomorphic to $\g$
itself occurs twice: once in the symmetric part ${\rm Sym}^2
(\g)$ of $\g^{\ot 2}$ and once in the skewsymmetric part
$\wedge^2(\g)$ (as usual, we assume that $\g$ acts onto itself
by the adjoint action and this action is extended to $\g^{\ot
2}$ via the Leibniz rule).  Denote these components by $\g_s$
and $\g_a$ respectively:
$$
\left.
\begin{array}{l}
\g_s\subset {\rm Sym}^2 (\g)\\
\g_a\subset \wedge^2(\g)
\end{array}\right\},\qquad \g_{a,s}\sim \g\;.
$$

Let us consider a non-trivial $\g$-morphism sending the
component $\g_a$ to $\g_s$ (it is unique up to a factor). Let us
extend this map to other components of $\wedge^2(\g)$ by 0. Let
$\{f,\,g\}_{inv}$ be the extension of this map  from
$\wedge^2(k(\g^*))$ to $k(\g^*)$ via the Leibniz rule. Then
there exist two values of $a$ such that the sum
\beq
\{f,\,g\}=\{f,\, g\}_r\,+\,a\{f,\,g\}_{inv}
\label{skobka}
\eeq
is a Poisson bracket (cf. \cite{DGS}). One of this two Poisson
brackets is a quasiclassical counterpart of the RE algebra (the
other one corresponds to a modified form of the RE algebra).  In
the sequel the appropriate $a$ is assumed to be fixed.  By this
the corresponding Poisson bracket is defined on $sl(n)^*$ 
(note that it is quadratic). In
order to pass to a Poisson bracket defined on $gl(n)^*$ we
should add one more generator, which Poisson commutes with all
other generators.  Let us observe that the bracket \r{skobka}
(or its extension to $gl(n)^*$) is compatible with the
corresponding linear Poisson-Lie bracket. The Poisson pencil
generated by these two brackets on $gl(n)^*$  is just the
quasiclassical counterpart of the  two parameter family $\Lqh$.

As for other simple Lie algebras $\g$  any
invariant correction to the bracket $\{\,\,,\,\,\}_r$ converting
it into  a quadratic  Poisson bracket (cf. \cite{DGS}) does not exist.

Now let us pass to the quasiclassical counterparts of the
quantum orbits above.  Observe that the Poisson bracket \r{skobka}
can be restricted to any orbit in $sl(n)^*$ (or in $gl(n)^*$ if
we take its extension), cf. \cite{D2} for a proof. 
 In particular, it is so for semisimple
generic orbits. Moreover, this restricted bracket is compatible
with the Kirillov-Kostant-Souriau one. The  Poisson pencil
generated by these two brackets is just the quasiclassical
counterpart of the algebra $k(M_{\mu}^{q\h})$ above.

Let us remark that for symmetric orbits the reduced Sklyanin
bracket becomes a particular case of this Poisson pencil. Thus,
for such type orbits the quantum objects can be described in two
ways: in terms of the Hopf-Galois extension or as a quotient of
the RE algebra or its quadratic-linear counterpart $\Lqh$.
However, the latter  way is more explicit and leads to objects
of "quantum affine algebraic geometry" (cf. \cite{DGK} where the
orbits of ${\bf CP}^n$ type were quantized in the spirit of such
a type geometry by means of an operator method).  For
non-symmetric orbits the quotients of the algebra $\Lq$ (or
$\Lqh$) and algebras arising from the reduced Sklyanin bracket
are completely different in spite of the fact that they both are
$\uq$-covariant.  Also note that the family of  $\uq$-covariant
algebras which are deformations of the coordinate ring of a
semisimple orbit in $gl(n)^*$ (or $sl(n)^*$) is large enough.
The reduced Sklyanin  bracket or the Poisson brackets
corresponding to the algebra $\Lq$ or $\Lqh$ represent only
particular cases of Poisson structures corresponding to this
family (cf. \cite{DGS} where such  Poisson structures are
classified).

\section{Basic line bundles over quantum orbits}
\setcounter{equation}{0}
In this section we introduce quantum line bundles over the
quantum orbits in question associated to the fundamental vector
$gl(n)$- (or what is the same $sl(n)$-) module $V$.  Let us fix
an  "orbit"  which is a quotient of one of the algebras $\Lq$,
$\Lqh$ or $\Lhh=U(gl(n)_{\h})$  over the ideal generated by
\r{center} (in the case of the  algebras $\Lq$, $\Lqh$ both in the
quasiclassical or non-quasiclassical cases).  

Similarly to the
above consideration we assume the roots $\mu_1,...,\mu_n$ of the
corresponding polynomial \r{poly} to be pairwise distinct.
In the sequel we use the notation $\kMm$ for such an orbit
(of one of the types above). Let
$$
\RVm=V\ot \kMm
$$
be a free right $\kMm$-module and be $\Rn$ its submodule
generated by the r.h.s. of \r{qb}. Let us consider the quotient
module $\RVm/\Rn$.
\begin{theorem}
The $\kMm$-module $\RVm/\Ru$ is non-trivial iff $\nu$ in \r{qb}
coincides with one of $\mu_i$, that is iff $\nu=\mu_i$ for some
$i$. In this case the module $\RVm/\Ru$ is projective. More
precisely,
\beq
\RVm/\Rmi = {\rm Im}\,P_i\;,\label{image}
\eeq
where
\beq
P_i = \prod_{j\not=i}^p\frac{L-\mu_j\,{\rm id}}{\mu_i-\mu_j}\;.
\label{L-proj}
\eeq
is a projector acting on the free module $\RVm$.  Furthermore,
in the quasiclassical case the module $\RVm/\Rmi$ is a flat
deformation of its classical counterpart. (In the case of the
algebra $U(gl(n)_{\h})$ we assume $L=A$.)
\end{theorem}
\begin{remark}
Let us note that we do not consider the $\kMm$-module
$\RVm/\Rmi$ as a quantum variety since the corresponding
"coordinate ring"  is not well defined. In order to introduce
such a ring we should define a commutation rule between the
space $V$ and the algebra $\kMm$.  However, apparently there is
no reasonable way to do it (if we want to preserve the flatness
of the deformation). We are planing to return to this question
in a future publication.
\end{remark}
{\bf Proof}\hspace*{5mm}
The proof of the theorem is based on the CH identities
\r{RE-CH}, \r{Lqh-CH} and \r{Ugl-CH}   and looks like that in the
classical case because the main difficulty is hidden in the
quantum version of the CH identity. The necessity of restriction
$\nu = \mu_i$ for some $i$ follows from relation \r{qb} and
${\overline P}(L)=0$ where ${\overline P}$ is defined by
\r{poly}. Note, that the latter relation can be rewritten in the
form:
$$
\prod_{i=1}^p(L-\mu_i \id)=0, \qquad p=\rank(R)\;.
$$

Now, by virtue of \r{qb} we have
$$
(\nu -\mu_1)\dots(\nu-\mu_p) v=0,\qquad v=(v_1,...,v_n)\;
$$
and if $\forall i \ \nu\not=\mu_i$ we have  $v=0$, that is the
module $\RVm/\Ru$ is trivial.  In order to prove the
non-triviality and projectivity of the module in the case
$\nu=\mu_i$ we consider the projectors \r{L-proj}.  Note that
the action of the projectors $P_i$ on the $\kMm$-module $\RVm$
is given by that of the matrix $L$ which is defined as follows
\beq
\sum_{i=1}^nv_i g^i(l)\triangleleft L=\sum_{i,j=1}^nv_j l^j_ig^i(l).
\label{modif}
\eeq
Thus, relation \r{qb} can be represented in the form
\beq
v\triangleleft L-\nu v,\; v\in V .
\label{qqb}
\eeq

Let us remark that the action \r{modif} of the matrix $L$ on the
space $V$ is coordinated with the coaction of the Hopf algebra
$\Tq$ in general and therefore (in the quasiclassical case) with
the action of the dual object, namely the QG $\uqs$. Taking
into account that $\mu_i$ are distinct we have the following
\begin{proposition}
The operators (\ref{L-proj}) form the full set of orthonormal
projective operators on the space $\RVm$, that is the following
properties hold
\begin{enumerate}
\item[{\bf i)}] $P_iP_j = \delta _{ij}P_i$,
\item[{\bf ii)}]$\displaystyle \sum_{i=1}^n P_i ={\rm id}$.
\end{enumerate}
\end{proposition}
{\bf Proof} is left to reader as an easy exercise.
\medskip
Let us return to the proof of the theorem. As follows from the
proposition the quotient $\RVm/\Rmi$ can be identified with $\Im
P_i$. This shows that the $\kMm$-module $\RVm/\Ru$ is
projective. Moreover, in the quasiclassical case it also implies
that this module is a flat deformation of its classical
counterpart since under a deformation the projectors are
deformed smoothly.  This completes the proof.

\section{Derived line bundles}
\setcounter{equation}{0}
In this section we consider the problem of constructing the
quantum line bundles different from the basic ones. We call them
{\em derived line bundles}.  First, consider the classical case.
Let us fix a generic semisimple orbit $M_{\mu}$ and  two line
bundles $E_{\mu_i},\,\, i=1,\,2$ (see Introduction).  Let us
consider their tensor product. We want to represent its
coordinate ring as $k(M_{\mu})$-module as well. It can be done
as follows.  

Let again $V$ be the fundamental vector
$sl(n)$-module. Consider the free right $k(M_{\mu})$-module
$$
\RVVM=\VV\ot k(M_{\mu})
$$
and extend the action of the matrix $L=A$ to the space $\VV$ by
setting
$$
(u\ot v)\triangleleft L=u\ot(v \triangleleft L)+(u \triangleleft L)\ot
v\qquad u, v \in V
$$
(some sort of the Leibniz rule). Thus, the extended matrix $L$
which will be denoted $L^{(2)}$ can be written as
$$
L_1+L_2,\quad L_1=\Id\ot L, \quad  L_2=P_{12}L_1P_{12}\;.
$$
Let us consider  the submodule $\Rn$ defined as in \r{qqb} but
with $v \in V$ replaced by $u\ot v$. Then the same problem
arises: for what value of $\nu$ the factor
$$
\RVVM/\Rn
$$
is not trivial.  It is not difficult to see that it is so iff
$\nu=\mu_i+\mu_j$ where $\mu_i$ are the roots of the polynomial
${\overline P}$, i.e., "eingenvalues" of the matrix $L$
corresponding to the orbit in question.  This allows us to find
the CH identity for the matrix $L^{(2)}$.

Also it can be found directly from the identity for the matrix
$L$ in the following way.  Let $\overline P(L)=0$  be the CH
identity for the matrix $L$ (in this section the coefficients of
polynomial $\overline P$ are supposed to be numerical). Then it
is not difficult to find an analogous relation for the matrix
$L^{(2)}$. For this it is sufficient to raise this matrix to the
powers $1,2,...,l(l+1)/2 $ where $l$ is the degree of the
polynomial $\overline P$ and on expressing the powers $L_1^l,
L_1^{l+1},...$ through $L_1^1,..., L_1^{l-1}$ and similarly for
the matrix $L_2$ we get the CH identity for the matrix
$L^{(2)}$.  The crucial property used in the construction is the
mutual commutativity of the matrices $L_1$ and $L_2$:
$$
L_1\,L_2= L_2\, L_1.
$$

Moreover, by assuming the sums $\mu_i+\mu_j$ to be pairwise
distinct we can construct the projector analogous to \r{proj}.
Nevertheless, we should restrict ourselves to the symmetric part
of the space $V^{\ot 2}$ to eliminate the multiplicity of the
quantity $\mu_i+\mu_j$ since (if $i\not=j$) it occurs once in
 the symmetric part of this space and once in its skewsymmetric
part. Thus, the projector corresponding to the eigenvalue
$\mu_i+\mu_j$ is the product of the projector onto symmetric
part of $V^{\ot 2}$ and that looking like \r{proj}.
Let us point out  that in a similar way it is possible to extend
the matrix $L$ to the higher tensor powers of the space $V$:
$V^{\ot l}, l=3, 4,...$. Thus, if $l=3$ the extended matrix is
defined as $L_1+L_2+L_3$ with the obvious definition of the
matrix $L_3$. The details are left to the reader.

Turn now to the quantum case, i.e., assume that the matrix $L$
is subject to relations \r{REA}. If we considered the matrix
$L^{(2)}=L_1+L_2$ with $L_2$ defined as above but with $P_{12}$
replaced by $R_{12}$ we would be unable to find the CH identity
for such an extension of the matrix $L$ since the matrix $L_2$
does not satisfy the polynomial relation valid for $L_1=L$. The
point is that the matrices $L_1$ and $L_2$ are not similar.
(However, if we chose as $L_2$ the matrix $R_{12}L_1R_{12}^{-1}$
then the matrices $L_1$ and $L_2$ would become similar but the
commutativity $L_1L_2=L_2L_1$ valid in the previous case by
virtue of the RE would be lost.) 

Nevertheless, we are interested
in an extension of the matrix $L$ to the symmetric part of the
space $V^{\ot 2}$. Let us define such an extension as follows
\beq
L_+=P_+L P_+\quad {\rm where}\quad P_{+}=\frac{q^{-1}\Id
+ R_{12}}{q+q^{-1}}\;.
\label{ext}
\eeq
Such a way to extend the matrix $L$ to the symmetric part of
$V^{\ot 2}$ is motivated by the following observation. In the
classical case $(q=1)$ such an extension of the matrix $L=A$
coincides (up to a factor, which does not matter for us) with
the restriction of the matrix $L_1+L_2$ to the symmetric part of
$V^{\ot 2}$.

In the sequel we will restrict ourselves to the case $\rank(R)=
2$.  This implies that the CH identity for  the matrix $L$ is
quadratic:
\beq
L^2-aL+b\,\id =0\qquad a=\mu_1+\mu_2,\ b=\mu_1\mu_2\;. \label{quadr}
\eeq
\begin{proposition}
If the CH identity for the matrix $L$ is of the form \r{quadr}
then the matrix $L_+$  defined by \r{ext} obeys the CH identity
of the form:
\bea
L_{+}^3-a(1+\frac{q^{- 1}}{2_q})\,L_{+}^2 + (a^2\frac{q^{-1}}{2_q}
- b)\,L_{+}+ab\frac{q^{-1}}{2_q}\,\id =0.
\eea
\end{proposition}
{\bf Proof}  \hspace*{5mm}
Taking into account the formulas for the symmetrizer $P_+^{(2)}$ we
can express the $R$-matrix via $P_+$  and rewrite the RE algebra in
the equivalent form:
\beq
P_{+}LP_{+}L - LP_{+}LP_{+} +\frac{q^{-1}}{2_q}(L^2P_{+} -
P_{+}L^2) =0\;. \label{RE-pr}
\eeq
If the matrix $L$ obeys the CH identity \r{quadr} we then have
\beq
P_{+}LP_{+}L - LP_{+}LP_{+} +a \frac{q^{-1}}{2_q}(LP_{+} -
P_{+}L) =0\;. \label{RE-ch}
\eeq
Now the CH identity for $L_+$ is a consequence of direct
calculations. Indeed, let us calculate successively the powers
of matrix $L_+$. Below we use the abbreviation $\xi =
a\frac{q^{-1}}{2_q}$.  For $L^2_+$ we get
$$
L_+^2\equiv ({P_+LP_+L})P_+ =({\rm use\  (\ref{RE-ch})}) =
LL_+ -\xi LP_+ +\xi L_+\;.
$$

We have here unwanted terms $LL_+$ and $LP_+$ and therefore
should calculate the next power of $L_+$ in order to get rid of
them. So
\begin{eqnarray*}
L_+^3&=&L_+^2L_+ = ({\rm insert\  L_+^2  \ above }) =LL_+^2
-\xi LL_+ +\xi L_+^2 =\\
&& ({\rm insert\  L_+^2 \ again\ and \ use \ (\ref{quadr})}) = \\
&& a (1+\frac{q^{-1}}{2_q}) LL_+ - a (1+\frac{q^{-1}}{2_q})
\xi LP_+ +(\xi^2 - b)L_+ +b\xi P_+.
\end{eqnarray*}
No new unwanted term besides $LP_+$ and $LL_+$ appears and
excluding them from expressions for $L_+^2$ and $L_+^3$ we come
to the following
$$
L_+^3 - a (1+\frac{q^{-1}}{2_q})L_+^2 +(a^2\frac{q^{-1}}{2_q} +b)
L_+ -ab\frac{q^{-1}}{2_q}P_+ = 0.
$$

It remains to observe that $P_+=\id$  on the symmetric part of
$V^{\ot 2}$.  In a similar way we can extend the matrix $L$ to
the higher symmetric powers of the space $V$: it suffices to
replace the projectors $P_+$ in the formula \r{ext} by the
symmetrizer in a given power.  However, the problem of finding
the corresponding CH identity is much more complicated and is
still open.

\section*{Appendix}
\def\theequation{A.\arabic{equation}}
\setcounter{equation}{0}
In this appendix we present the explicit form of the
coefficients entering  (\ref{Lqh-CH}) and show that in the
quasiclassical case these coefficients have finite limit as
$q\to 1$. This fact gives rise to the formula (\ref{Ugl-CH}).
Let us recall that  in general $p={\rm rank}(R) \not=n=\dim V$.
Also recall that the degree of the polynomial $\overline P$  and
hence the number of the roots $\mu_i$ is equal to $p$.  We will
use the fact, that the $\Lqh$ algebra (\ref{Lqh-com}) can be
formally obtained from $\Lq$ by the shift of generators $l^i_j =
{\bar l}^i_j - \delta^i_j h$ and subsequent changing of the
parameter: $\h=h\lambda$. So, if we realize these two operations
in the CH identity (\ref{RE-CH}) for the matrix $L$ we get some
polynomial identity on the matrix ${\bar L}$ formed by the
generators of the algebra $\Lqh$.  

The problem is to find in
which way the central elements $\sigma_k(L)$ are transformed.
This can be directly calculated from the definition of
$\sigma_k(L)$ (\ref{L-sig}). Indeed, on making the above
mentioned shift in $\sigma_k(L)$ we have to transform the number
of arising terms, their typical form being
\beq
{\rm Tr}_{(12\dots p)}P_-^{(p)}(L_1R_{1}\dots R_{k-1})^sR_1\dots
R_{k-1}(L_1R_{1}\dots R_{k-1})^{k-s-1}\;.
\label{cancel}
\eeq
In the above formula the concise notation $R_i\equiv R_{i\,i+1}$
is used.  Now one should "draw out" the string of $R$-matrices
and cancel them on the projector $P_-^{(p)}$. Then it is
necessary to get rid of all the matrices $R_{k-1}$ in monomials
$(L_1R_{1}\dots R_{k-1})$.  The basic formulas for such
transformations are
$$
({\bar L}_1R_{1}\dots R_{k})R_{i} =
R_{i+1} ({\bar L}_1R_{1}\dots R_{k})\qquad \forall\,i\le k-1\, ,
$$
(this is a trivial consequence of the Yang-Baxter equation) and
$$
({\bar L}_1R_{1}\dots R_{k})^r = ({\bar L}_1R_{1}\dots
R_{k-1})^rR_kR_{k-1}\dots  R_{k-r+1}
$$
which in turn is a direct consequence of the previous relation.
The cancellation of $R$-matrices in \r{cancel} is due to the
cyclic property of trace and the defining property of the
skewsymmetrizer $P_-^{(r)}$:
$$
R_i\,P_-^{(r)} = P_-^{(r)}\,R_i = -\frac{1}{q}\, P_-^{(r)}\quad
\forall i\le r-1\;.
$$

Now after straightforward calculations we come to the following
transformation of coefficients:
\beq
\sigma_k(L) \longrightarrow \sum_{r=0}^k \left(-\frac{\hbar}{\lambda}
\right)^r q^{-r(p-1)}C^p_k\,\frac{[C^k_r]_q}{[C^{k-r}_p]_q}\,
\sigma_{k-r}({\bar L})\;.
\label{sig-izm}
\eeq
Here the symbol $[C^k_p]_q$ stands for the $q$-binomial
coefficient
$$
[C^k_p]_q = \frac{p_q!}{k_q!(p-k)_q!}
$$
where $p_q! = (p-1)_q!p_q$ and $q$-numbers $r_q$ are defined as
$$
r_q\equiv \frac{q^r - q^{-r}}{q-q^{-1}}\;.
$$

The elements $\sigma_k({\bar L})$ are defined in (\ref{L-sig})
where matrix $L$ should be changed for ${\bar L}$. Obviously,
these elements are central in $\Lqh$ algebra.  Now given the
rule (\ref{sig-izm}) it is not so difficult to obtain from
(\ref{RE-CH}) the CH identity for the $\Lqh$ algebra
$$
(-{\bar L})^p + \sum_{k=0}^{p-1} (-{\bar L})^k
\sigma^{(\hbar)}_{p-k}({\bar L}) \equiv
0\;.\eqno{(\ref{Lqh-CH})}
$$
The quantities $\sigma^{(\hbar)}_{p-k}({\bar L})$ are the
following polynomials in $\h$ with coefficients depending on
$\sigma_{p-k}({\bar L})$:
\beq
\sigma^{(\hbar)}_{p-k}({\bar L}) = \sigma_{p-k}({\bar L}) +
\sum_{r=1}^{p-k}  \hbar^r\omega^{(p)}_{r+k,k}\,
\sigma_{p-k-r}({\bar L})\;.
\label{sig-h}
\eeq
The numeric coefficients $\omega^{(p)}_{s,k}$ are as follows  ($s>k$)
\beq
\omega^{(p)}_{s,k} =
\frac{\lambda^{k-s}}{[C^s_p]_q}\,\sum_{r=0}^{s-k}
(-1)^rq^{-r(p-1)} C^k_{s-r}C^r_{p-s+r}
[C^{s-r}_p]_q\;.\label{om-def}
\eeq

Now we calculate the sum in (\ref{om-def}) and show that it is
proportional to $\lambda^{s-k}$  and therefore the whole
coefficient $\h^r\omega^{(p)}_{r-k,k}$ admits a non-singular classical
limit as $q\rightarrow 1,\,\,\h=$const. On taking this limit in CH
identity
(\ref{Lqh-CH}) we get the corresponding identity (\ref{Ugl-CH})
for the matrix $A$ of $\Lh$ algebra (\ref{U-com}).  Denote the
sum in (\ref{om-def}) by $\xi^{(p)}_{s,k}$:
$$
\xi^{(p)}_{s,k}\equiv \sum_{r=0}^{s-k}
(-1)^rq^{-r(p-1)} C^k_{s-r}C^r_{p-s+r} [C^{s-r}_p]_q
$$
and calculate the generating function $\Phi_{\xi}(x,y)$ of the
coefficients $\xi^{(p)}_{s,k}$
$$
\Phi_{\xi}(x,y)\stackrel{\mbox{\tiny def}}{=}
\sum_{s=0}^p\sum_{k=0}^s (-x)^{p-s}\,(-y)^k\,\xi^{(p)}_{s,k}\;.
\label{def-Phi}
$$

If one knew the function $\Phi_{\xi}(x,y)$ then the coefficients
$\xi^{(p)}_{s,k}$ could be found as
\beq
\xi^{(p)}_{s,k} = \frac{(-1)^{p-s+k}}{k!(p-s)!}\,
\left[\frac{\partial^{p-s}}{\partial x^{p-s}}
\frac{\partial^{k}}{\partial y^{k}}\,
\Phi_{\xi}(x,y)\right]_{x=y=0} \;.
\label{ksi-d-Phi}
\eeq
The calculation of $\Phi_{\xi}(x,y)$ is rather simple where the
only thing we need for is the Newton binomial formula and its
$q$-analogue
$$
\sum_{k=0}^p x^k q^{-k(p-1)}[C^k_p]_q =
\prod_{k=0}^{p-1}\left( 
1 + \frac{x}{q^{2k}}\right) \;. 
$$ 

So we present the final result  
\beq 
\Phi_{\xi}(x,y) = (-1)^pq^{-p(p-1)}\prod_{k=0}^p
(xq^{p-1}+yq^{2k} - \lambda\,q^kk_q)
\label{Phi-rez}\;. 
\eeq  
At last, upon taking the partial derivatives in
(\ref{ksi-d-Phi}) we find the form of coefficients
$\xi^{(p)}_{s,k}$:
\beq 
\xi^{(p)}_{s,k} = \lambda^{s-k}q^{{(p-1)(p-2s)\over 2}}
(p-1)_q!\, (V_k+V_s)\;, 
\eeq 
where by $V_{k}$ and $V_{s}$ we denote the following sums 
$$ 
V_k\equiv \left\{ 
\begin{array}{l} 
\qquad 0\quad {\rm if}\ k=0\\ 
 \\ 
\displaystyle 
\sum_{1\le l_1<\dots<l_{k-1}\le p-1 \atop  
{1\le r_1<\dots<r_{p-s}\le p-1 \atop  \{l_i\}\cap \{r_j\} = 
\emptyset}} \frac{q^{(l_1+\dots + l_{k-1} - r_1 - \dots  
- r_{p-s})}}{(l_1)_q\dots (l_{k-1})_q(r_1)_q\dots (r_{p-s})_q} 
\end{array}\right.\quad {\rm if}\ k\not=0 
$$ 
$$ 
V_s\equiv \left\{ 
\begin{array}{l} 
\qquad 0
\quad {\rm if}\ s=p\\ 
 \\ 
\displaystyle 
\sum_{1\le l_1<\dots<l_{k}\le p-1 \atop  
{1\le r_1<\dots<r_{p-s-1}\le p-1 \atop  \{l_i\}\cap \{r_j\} = 
\emptyset}} \frac{q^{(l_1+\dots + l_{k} - r_1 - \dots  
- r_{p-s-1})}}{(l_1)_q\dots (l_{k})_q(r_1)_q\dots (r_{p-s-1})_q} 
\end{array}\right.\quad {\rm if}\ s\not=p\;. 
$$  
In particular, 
$$ 
\xi^{(p)}_{s,s}\equiv 1\qquad \xi^{(p)}_{p,0}\equiv 0\;. 
$$ 

Now we can find the coefficients $\omega^{(p)}_{s,k}$ entering the 
CH identity (\ref{Lqh-CH}) of the $\Lqh$ algebra 
$$ 
\omega^{(p)}_{s,k} \equiv 
\frac{\lambda^{k-s}}{[C_p^s]_q}\xi^{(p)}_{s,k} 
=  q^{{(p-1)(p-2s)\over 2}}\frac{(p-1)_q!}{[C_p^s]_q}\, 
(V_k+V_s)\;. 
$$ 
This gives the explicit form of $\sigma^{(\hbar)}_k$ and
completes the proof of the CH identity for the algebra $\Lqh$.
From the above relation it is obvious that the coefficients
$\omega^{(p)}_{s,k}$  admits a non-singular classical limit
\beq 
\lim_{q\rightarrow 1} \omega^{(p)}_{s,k} \equiv \rho^{(p)}_{s,k} 
= \frac{(p-1)!}{C_p^s}\, (V_k^{cl}+V_s^{cl})\;,
\label{ro-def} 
\eeq 
where $V_{k,s}^{cl}$ are given by formulas for $V_{k,s}$ with
substitution $q=1$ and all $q$-numbers changed for ordinary
ones.  

Finally, on taking into account that $\Lh =
\lim_{q\rightarrow 1}\Lqh$ we deduce from (\ref{Lqh-CH}) the CH
identity \r{Ugl-CH} where coefficients $\tau^{(\hbar)}_{p-k}(A)$
are the classical limit of $\sigma^{(\hbar)}_{p-k}({\bar L})$ in
(\ref{sig-h}):
\beq 
\tau^{(\hbar)}_{p-k}(A) = \lim_{q\rightarrow 1} 
\sigma^{(\hbar)}_{p-k}({\bar L}) = \tau_{p-k}(A) + 
\sum_{s=1}^{p-k}  \hbar^s\rho^{(p)}_{s+k,k}
\, \tau_{p-k-s}(A)\label{A9}
\eeq
with $\rho^{(p)}_{s,k}$ defined in (\ref{ro-def}). The central
elements $\tau_k(A)$ have the form:
$$
\tau_k(A) = \frac{C^k_p}{p!}\, \varepsilon_{i_1\dots i_k
a_{k+1}\dots a_p} A^{i_1}_{j_1}\dots A^{i_k}_{j_k}
\varepsilon^{j_1\dots j_k a_{k+1}\dots a_p}\;,
$$
$\varepsilon$ being the skew-symmetric Levi-Civita tensor.
These elements are analogues of spectral invariants of the usual
matrix with commutative entries: in that case each $\tau_k$ is
the sum of all principal minors of $k$-th order.  

From identity
(\ref{A9}) one can see, that in the classical
non-commutative case each coefficient in the CH identity is
modified by adding a polynomial in $\hbar$. In particular, the
free term of the identity which can be treated as
non-commutative analogue of the determinant is given by the
following formula
\beq
\det A +\sum_{k=1}^{p-1}\hbar^k \rho^{(p)}_{k,0} \tau_{p-k}(A).
\eeq
In this formula it is taken into account that
$\rho^{(p)}_{p,0}=0$.


\begin{thebibliography}{DDD}

\bibitem[AG]{AG}
Akueson P. and Gurevich D. {\it Some aspects of braided geometry:
differential calculus, tangent space, gauge theory}, J. Phys. A:
Math Gen. {\bf 32} (1999), pp 4183--4197.

\bibitem[BL]{BiLo}
Biedenharn L.C. and Louck J.D. {\em A pattern calculas for tensor
operqtors in the unitary groups}, CMP {\bf 8}  (1968), pp 89--131.\par
Biedenharn L.C. and Louck J.D. {\it Canonical unit adjoint tensor
operators in $U(n)$}  J. Math. Phys. {\bf 11}  (1970),
pp 2368-2414.

\bibitem[D1]{D1}
Donin J. {\it  Double quantization on the coadjoint representation
of $sl(n)$},  Czech. J. Phys. {\bf 47} (1997), pp 1115--1122.

\bibitem[D2]{D2} Donin J. {\em Double quantization on the coadjoint
representations of simple Lie groups and their orbits}, Preprint
MPIM-Bonn
 99-103.

\bibitem[DG]{DG}
Donin J. and Gurevich D. {\it Some Poisson structures associated
to Drinfeld-Jimbo R-matrices and their quantization}, Israel
Math. J. {\bf 92} (1995), pp 23--32.

\bibitem[DGK]{DGK}
Donin J., Gurevich D. and Khoroshkin S. {\it Double quantization of
${\bf CP}^n$ type orbits by generalized Verma modules},
Jour. Geom. Phys. {\bf 28} (1998), pp 384--406.

\bibitem[DGS]{DGS} Donin J., Gurevich D. and Shnider S.
{\em Double quantization
 on some orbits in the coadjoint representations of simple Lie
 groups},
CMP {\bf 204} (1999), pp 39--60.

\bibitem[FRT]{FRT}
Faddeev L.D., Reshetikhin N.Yu and Takhtajan L.A.
{\it Quantization of Lie groups and Lie algebras},
 Leningrad Math. J. {\bf 1} (1990), pp 193--226.

\bibitem[G-T]{GKLLRT}
Gelfand I.M., Krob D., Lascoux A., Leclerc
B.,  Retakh V.S. and Thibon J.-Y. {\it Noncommutative symmetric
functions}, Adv. Math. {\bf 112} (1995), pp 218--348.

\bibitem[G]{G}
Gurevich D. {\it Algebraic aspects of the quantum Yang-Baxter
equation}, Leningrad Math. J. {\bf 2} (1991), pp 801--828.

\bibitem[GR]{GR}
Gurevich D. and Rubtsov V. {\it Quantization of Poisson pencils
and generalized Lie algebras}, Teor.  Mat. Phys. {\bf 103}
(1995), pp 476--488.

\bibitem[GPS]{GPS}
Gurevich D.I., Pyatov P.N. and Saponov P.A.
{\it Hecke Symmetries and Characteristic Relations on Reflection
Equation Algebras},  Lett. Math. Phys. {\bf 41} (1997), pp 255--264.

\bibitem[HM]{HM}
Hajac P. and Majid S. {\em Projective module description of the
Q-monopole},  CMP {\bf 206} (1999), pp 247-264.

\bibitem[IOP]{IOP}
Isaev A., Ogievetsky O. and Pyatov P. {\it Generalized
Cayley-Hamilton-Newton identities}, Czech. J. Phys. {\bf 48}
(1998), pp 1369--1374.

\bibitem[KSk]{KulSkl}
Kulish P.P. and Sklyanin E.K. {\it Algebraic structures related to
reflection
 equation}, J. Phys. A: Math. Gen. {\bf 25}
(1992), pp 5663--5975.

\bibitem[KSa]{KulSas}
Kulish P.P. and Sasaki R. {\it Covariance properties of Reflection
equation algebras}, Prog. Theor. Phys. {\bf 89} (1993), pp 741--761.

\bibitem[M1]{M1}
Majid Sh. {\it Examples of braided groups and braided matrices},
J. Math. Phys. {\bf 32} (1991), pp 3246--3252.

\bibitem[M2]{M2}
Majid Sh. {\it Foundations of Quantum Group theory},
Cambridge University Press 1995.

\bibitem[Sh]{Sh}
Schneider S. {\it Hopf Galois Extensions, Crossed Products,
and Clifford Theory} in: Bergen. J., Montgomery S. (eds.) Advances in
Hopf Algebras. Lecture Notes in Pure and Applied Mathematics. Marcel
Dekker Inc. {\bf 158} (1994), pp 267--297.

\bibitem[Se]{Se}
Serre J.-P.
{\em Modules projectifs et espaces
fibr\'es a fibre vectorielle},
 S\'eminaire Dubreil-Pisot
Fasc. 2, Expos\'e 23 (1957/1958). 

\bibitem[Sw]{Sw}
Swan R. {\em Vector bundles and projective modules}, Trans. Am. Math.
Soc. {\bf 105} (1962), pp 264--277.
\end{thebibliography}
\end{document}